# An improved lower-bound on the counterfeit coins problem

Li An-Ping
Beijing 100085, P.R. China
apli0001@sina.com

Hagen von Eitzen
hagen@von-eitzen.de

Abstract

In this paper, we will give an improvement on the lower bound for the counterfeit coins problem in the case that the number of false coins is unknown in advance.

Searching for counterfeit coins in a number of coins with same semblance by a balance is a well-known combinatorial search problem called the counterfeit coins problem. The problem has a longer history and gotten a intensive researches, for the detail refer to see the papers [1]~[11]. The problem has several versions, e.g. the number of the fakes is assumed known in advance, and/or the fakes are known lighter or heavier than the normals. Usually, it is assumed that the coins will be permitted to be remarked by numbers in order to be distinguished each other.

In paper [5], we discussed a general case that the number of the fakes is unknown beforehand.

Suppose that $\mathcal{S}$ is a set of $n$ coins with same semblance, in which possibly there are some counterfeit coins, which are heavier (or lighter) than the normals. Denoted by $g(n)$ the least number of weighings need to find all the fakes in $\mathcal{S}$ by a balance, assumed that additional normal coins will be available if needed. Our main result in [5] is as following

$$\lceil n \cdot \log_3 2 \rceil \leq g(n) \leq \lceil 7n/11 \rceil, \tag{1}$$

but with a exception that $g(3) = 3$.

Clearly, the lower-bound, i.e. the left-hand side of (1) is just the information theoretic bound of $g(n)$, in this paper, we will give an improvement over it, our main result is that

**Proposition 1**
$$g(n) \geq \lceil \log_3(2^n + 2^{n-5} + 2^{n-6} + 2^{n-7} + 2^{n-13} + 2^{n-14} + 2^{n-16} + 2^{n-17}) \rceil. \tag{2}$$

At first, we introduce some notations. Let $\mathcal{A}$ be the set of the three symbols $<, =, >$, i.e, $\mathcal{A} = \{<, =, >\}$, a vector $\alpha = (v_1, \cdots, v_k) \in \mathcal{A}^k$ will be called a direction, and a subset $X \subseteq \mathcal{S}$ will be called a objective in the direction $\alpha$, if the state of the balance in the i-th weighing for $X$ is equal to $v_i, 1 \leq i \leq k$. Denoted by $\mathcal{S}_\alpha$ the set of all the objectives in the direction $\alpha$.

As usually, for a set $A$, $2^A$ and $\binom{A}{k}$ will represent the set of all subsets of $A$ and the set of all the $k$-subsets of set $A$ respectively. Similar to Cartesian product, for any two subsets $\Delta, \Gamma \subseteq 2^{\mathcal{S}}$, we define $\Delta \cdot \Gamma = \{A \cup B \mid A \in \Delta, B \in \Gamma\}$.

In this paper, we will employ a combinatorial identity stated in the following lemma 1.

**Lemma 1.** Suppose that $s, t$ are two non-negative integers, and $k$ is an integer, then

$$\sum_i C_s^i \cdot C_t^{i+k} = C_{s+t}^{t-k}. \tag{3}$$

The formula (3) may be obtained by a basic combinatorial calculus, which has been omitted.

For a subset $\Gamma \subseteq \mathcal{S}$, denoted by $\mathcal{S}_\alpha(\Gamma) = \{X \mid X \in \mathcal{S}_\alpha, X \subseteq \Gamma\}$. Suppose that $L_i : R_i$, $i = 1, 2, \ldots, k$, are the first $k$ weighings in a direction $\alpha$, denoted by $\Gamma = \bigcup_{1 \leq i \leq k} (L_i \cup R_i)$, then it is easy to know that

$$|\mathcal{S}_\alpha| = |\mathcal{S}_\alpha(\Gamma)| \cdot 2^{n-|\Gamma|}. \tag{4}$$

For an integer $k$, define $\zeta_k = \max\{|\mathcal{S}_\alpha| \mid |\alpha| = k\}$.

**Lemma 2.**

$$\zeta_2 \geq 2^n \times 15/128. \tag{5}$$

Proof. Suppose that the first weighing is that $A : B$, and the second weighing is that $L : R$. Denoted by $\Gamma = A \cup B$, $C_1 = L \setminus \Gamma$, $C_2 = R \setminus \Gamma$, $|C_1| = s$, $|C_2| = t$. For an integer $i$, denoted by

$$\sigma_i(s,t) = \sum_k \binom{C_1}{k} \binom{C_2}{k+i}, \tag{6}$$

By Lemma 1, it has

$$|\sigma_i(s,t)| = \sum_k \binom{|C_1|}{k} \binom{|C_2|}{k+i} = C_{s+t}^{t-i} \tag{7}$$

We will also simply write $\sigma_i(s,t)$ as $\sigma_i$ in some apparent cases. Moreover, for a subset $X \subseteq \Gamma$, denoted by $\delta(X) = |X \cap L| - |X \cap R|$. Let $\Lambda = \Gamma \cup C_1 \cup C_2$, $v \in \{<, =, >\}$, it is clear that

$$\mathcal{S}_{(v,=)}(\Lambda) = \sum_{X \in \mathcal{S}_{(v)}(\Gamma)} \sigma_{\delta(X)} \cdot (X \cap L) \cdot (X \cap R). \tag{8}$$

So,

$$|\mathcal{S}_{(v,=)}(\Lambda)| = \sum_{X \in \mathcal{S}_{(v)}(\Gamma)} |\sigma_{\delta(X)}| = \sum_{X \in \mathcal{S}_{(v)}(\Gamma)} C_{s+t}^{t-\delta(X)} \leq |\mathcal{S}_{(v)}(\Gamma)| \times C_m^{[m/2]}, \tag{9}$$

where $m = s + t$.

It is easy to know that

$$\max\{|\mathcal{S}_{(<)}(\Gamma)|,|\mathcal{S}_{(=)}(\Gamma)|,|\mathcal{S}_{(>)}(\Gamma)|\} \geq \begin{cases} 1, & \text{if } |\Gamma|=1, \\ 2, & \text{if } |\Gamma|=2, \\ 4, & \text{if } |\Gamma|=3. \end{cases} \qquad (10)$$

Moreover, if $|\Gamma|=5$, then

$$\max\{|\mathcal{S}_{(<)}(\Gamma)|,|\mathcal{S}_{(>)}(\Gamma)|\} \geq 16. \qquad (11)$$

On the other hand, if $|\Gamma| \geq 6$, then

$$\max\{|\mathcal{S}_{(<)}(\Gamma)|,|\mathcal{S}_{(>)}(\Gamma)|\} \geq (2^{|\Gamma|} - C_{|\Gamma|}^{[|\Gamma|/2]})/2. \qquad (12)$$

Hence, the remain cases to be checked are

i) $|\Gamma|=4, |A|=|B|=2, m \leq 4$,
ii) $|\Gamma|=6, |A|=|B|=3, m \leq 4$.

For the case i), it has that

$$\max\{|\mathcal{S}_{(=,<)}(\Lambda)|,|\mathcal{S}_{(=,=)}(\Lambda)|,|\mathcal{S}_{(=,>)}(\Lambda)|\} \geq 2^{|A|}/8, \qquad (13)$$

For the case ii), it has that

$$\max\{|\mathcal{S}_{(<,<)}(\Lambda)|,|\mathcal{S}_{(<,=)}(\Lambda)|,|\mathcal{S}_{(<,>)}(\Lambda)|\} \geq 15 \times 2^{|A|}/128. \qquad (14)$$

The estimation (5) is from (10) to (14). □

**Lemma 3.**
$$\zeta_4 \geq 2^{n-18} \times 3414.$$

Proof. Follow with the proof of Lemma 2, suppose that the third weighing and the fourth are $U:V$ and $X:Y$ respectively. From the proof above, we know that

$|A|=|B|=3$ and $|(L \cup R) \setminus (A \cup B)| \leq 4$.

Similarly, there are

$|(U \cup V) \setminus (A \cup B \cup L \cup R)| \leq 4$ and $|(X \cup Y) \setminus (A \cup B \cup L \cup R \cup U \cup V)| \leq 4$.

Moreover, by Lemma 2, it has

$$\zeta_3 \geq 2^{n-14} \times 640.$$

And so

$$\zeta_4 \geq 2^{n-14-a} \times (2^a \times 640 + \varepsilon)/3,$$

$0 \leq a \leq 4$, $\varepsilon = 1$, or $2$, as $a$ is odd, or even. It follows

$$\zeta_4 \geq 2^{n-18} \times 3414. \qquad (16)$$

□

Proof of Proposition 1. Let $3^{k-1} < 2^n \leq 3^k$, it is known that the information theoretic bound of $g(n)$ is equal to $\lceil \log_3 2^n \rceil$, hence $g(n) \geq k$. Denoted by

$$\omega = 2^{n-5} + 2^{n-6} + 2^{n-7} + 2^{n-13} + 2^{n-14} + 2^{n-16} + 2^{n-17},$$

and suppose that $3^k = 2^n + x$, if $\omega \leq x$, then

$$2^n + \omega \leq 2^n + x = 3^k.$$

Namely,

$$g(n) \geq k = \lceil \log_3(2^n + \omega) \rceil. \tag{17}$$

Thereby, we assume that $x < \omega$, then $3^k < 2^n + \omega$, and

$$2^n > 3^k \times \frac{2^{17}}{2^{17} + 7195}. \tag{18}$$

By Lemma 3 and (18), it has

$$\zeta_4 \geq 2^n \times \frac{1707}{2^{17}} > 3^k \times \frac{2^{17}}{2^{17} + 7195} \times \frac{1707}{2^{17}} = 3^{k-4}.$$

This means that

$$g(n) \geq k + 1 = \lceil \log_3(2^n + \omega) \rceil. \tag{19}$$

The proof of Proposition 1 has been finished. □

*Remark.* With computer aid, it has shown that for the second weighing follow the direction $(<)$, there are following six types of ones which give $\zeta_2 = 2^{n-7} \times 15$ up to the equivalency,

1) $\{1,2,4\} : \{3,7,8\}$
2) $\{1,4,5\} : \{6,7,8\}$
3) $\{1,4,7\} : \{2,5,8\}$
4) $\{1,7,8\} : \{2,9,10\}$
5) $\{4,7,8\} : \{5,9,10\}$
6) $\{7,8,9\} : \{1,4,10\}$

With Lemma2 and a short program, it may be known that the second weighing of the first five types all give the estimation

$$\zeta_3 = 2^{n-10} \times 40, \tag{20}$$

and

$$\zeta_4 = 2^{n-14} \times (640+2)/3 = 2^{n-13} \times 107.  \qquad (21)$$

For the sixth type of the second weighing, follow the direction (<,<), with the third weighing $\{1,2,3,10,11,13\}$ : $\{5,6,7,8,9,12\}$ give the estimation $\zeta_3 = 2^{n-13} \times 320 \,(= 2^{n-10} \times 40)$, and then

$$\zeta_4 = 2^{n-17} \times 1712 = 2^{n-13} \times 107.  \qquad (22)$$

Overall, it follows

$$\zeta_3 = 2^{n-7} \times 5, \quad \zeta_4 = 2^{n-13} \times 107.  \qquad (23)$$

So,

$$g(n) \geq 4 + \lceil \log_3 \zeta_4 \rceil \geq \lceil n \cdot \log_3 2 + \log_3(107 \times 81/2^{13}) \rceil,$$

or,

$$g(n) \geq \lceil \log_3(2^n + 2^{n-5} + 2^{n-6} + 2^{n-7} + 2^{n-9} + 2^{n-10} + 2^{n-12} + 2^{n-13}) \rceil,  \qquad (24)$$

which is a little better than (2).